\documentclass{amsart}

\usepackage{amssymb, bm, amsfonts, stmaryrd, amsmath, amsthm, enumerate, mathtools, mathrsfs, colonequals, graphicx, color, nccmath, thmtools}

\usepackage[hidelinks]{hyperref}
\usepackage[capitalise]{cleveref}
\usepackage[parfill]{parskip}
\usepackage[margin=1.2in]{geometry}
\usepackage{multicol}
\usepackage{wrapfig,tikz,tikz-cd}
\usetikzlibrary{arrows}
\usepackage[foot]{amsaddr}

\usepackage[normalem]{ulem}

\newcommand{\rightcosets}[2]{%
  \mathord{\mathpalette\rightcosetslower{#1}\backslash #2}%
}
\newcommand{\rightcosetslower}[2]{%
  \raisebox{-0.3\height}{$#1#2$}%
}

\theoremstyle{plain}\newtheorem{theorem}{Theorem}[section]

\theoremstyle{plain}\newtheorem{conjecture}[theorem]{Conjecture}

\theoremstyle{plain}\newtheorem{corollary}[theorem]{Corollary}

\newtheorem*{corollary*}{Corollary}

\theoremstyle{plain}\newtheorem{lemma}[theorem]{Lemma}
\theoremstyle{plain}\newtheorem{proposition}[theorem]{Proposition}
\theoremstyle{plain}
\theoremstyle{definition}\newtheorem{definition}[theorem]{Definition}
\theoremstyle{definition}
\theoremstyle{definition}
\theoremstyle{definition}
\theoremstyle{definition}
\theoremstyle{definition}
\theoremstyle{definition}
\theoremstyle{definition}\newtheorem{propositiondefinition}[theorem]{Definition/Proposition}

\DeclareMathOperator{\HH}{HH}
\DeclareMathOperator{\Hom}{Hom}

\DeclareMathOperator{\Br}{Br}

\renewcommand{\phi}{\varphi}
\newcommand{\F}{\mathcal{F}}

\title{On Non-vanishing of the first Hochschild cohomology of blocks.}
\author{Matthew Antrobus$^1$}
\address{
$^1$Department of Mathematics, University of Manchester, Manchester M13 9PL, United Kingdom {\tt matthew.antrobus@manchester.ac.uk}}

\thanks{\noindent
{\it Funding.}
The author was funded by the Heilbronn Institute for Mathematical Research. 
}

\begin{document}
\begin{abstract}
	Let $B$ be a block of a finite group algebra $kG$, where $k$ is an algebraically closed field of characteristic $p$. Inspired by the results of \cite{FJL}, Linckelmann conjectured that $\HH^1(B)\not=0$ when $B$ has non-trivial defect. By developing the methods of \cite{BR} we give a fusion theoretic criterion for non-vanishing, proving the conjecture holds for all blocks of abelian defect. Further, we show the criterion holds for all blocks at $p>5$, proving the conjecture for all such primes.
\end{abstract}

\maketitle

\section{Introduction}



Let $G$ be a finite group, and $k$ be a field whose characteristic $p>0$ divides the order of $G$. The first Hochschild cohomology of the group algebra $kG$ was shown to be non-vanishing in \cite{FJL}, verifying a consequence of the Donald-Flanigan conjecture. Linckelmann conjectured the following strengthening of this result.

\begin{conjecture}[\cite{LinckelmannTalk}, \cite{Linckelmann2}]\label{markus1}
If $B$ is a block of a finite group algebra $kG$ with a non-trivial defect group, then $\HH^1(B)\not=0$.
\end{conjecture}

Our main theorem proves the conjecture in characteristic greater than $5$.

\begin{theorem}\label{nonvanishing}
    Let $G$ be a finite group, and $k$ be a field of characteristic greater than $5$. Take $B$ a block of $kG$ with a non-trivial defect group. Then $\HH^1(B)\not=0$. 
\end{theorem}

We are also able to verify the conjecture for blocks with non-trivial abelian defect at all primes.

\begin{restatable}{theorem}{AbelianTheorem}\label{AbelianTheorem} Suppose $B$ is a block of $kG$ with non-trivial abelian defect. Then $\HH^1(B)\not=0$.
\end{restatable}

The conjecture has been verified for various families of groups and blocks. In \cite[Theorem 1.2]{BKL2}, Benson, Kessar and Linckelmann calculated the dimensions of $\HH^1(B)$ for blocks of symmetric groups, and in \cite[Theorem 1.1]{Murphy2} Murphy found the dimensions corresponding to blocks of the Mathieu groups. Most blocks of the Janko groups were shown to have non-vanishing first Hochschild cohomology by Serwene and Todea in \cite[Theorem 1.6]{Todea2}, as well as some blocks of the alternating groups. All blocks of alternating groups, as well as blocks of groups of Lie type in defining characteristic, were shown to satisfy Conjecture \ref{markus1} in \cite[Theorem B]{BR} by Briggs, Rubio y Degrassi and the author. In \cite[Theorem A]{WLA}, Wen, Rubio y Degrassi and the author verified the conjecture for all blocks of sporadic simple groups, with exact dimensions given for most blocks.

Blocks with normal defect were shown to obey the conjecture by Todea in \cite[Corollary 1.4]{Todea1}, as were blocks of $p$-solvable groups, and inertial blocks. Moreover \cite[Theorem 1.5]{Todea2} demonstrated that should a counterexample obey certain strong properties, then there exists a counterexample amongst blocks of quasisimple groups. Blocks with maximal defect at $p>5$ were shown to have non-vanishing first Hochschild cohomology in \cite{BR}, utilising the classification of simple groups.

Let $B=kGb$ be a block of $kG$ with defect group $D$. As Hochschild cohomology respects extension of scalars, we may assume without loss of generality that $k$ is algebraically closed. The non-vanishing result of \cite{BR} relies on \emph{extra non-Schur elements} of $G$ lying in $D$, inspired by weaker properties in \cite{FJL}. Extending this definition to fusion systems, we define fusion non-Schur elements of a fusion system.

 \begin{restatable}{definition}{MainDef} Given a fusion system $\mathcal{F}$ on $P$, we say $x \in P$ is a fusion non-Schur element when $\langle x \rangle \leq P$ is fully-$\mathcal{F}$-centralized and $x \not \in \Phi(C_P(x))$. We write $\mathrm{NS}(\F)$ when such an element exists.
 \end{restatable} 
 
Extra non-Schur elements were utilised in \cite{BR} to transfer group theoretic information from $G$ to the first Hochschild cohomology of blocks related to $G$, using the Batalin-Vilkovisky (BV) operator $\Delta^1_{kG}:\HH^1(kG) \rightarrow \HH^0(kG)$ on Hochschild cohomology. Using techniques from block theory, they showed that the presence of an extra non-Schur element of $G$ lying in $D$ implies that $\operatorname{Im}(\Delta^1_{kG})$ has non-zero projection onto $\HH^0(B)$. When the defect group is maximal, the existence of an extra non-Schur element in $D$ can be reduced to simple groups. However, when the defect is not maximal, their methods are unable in general to demonstrate existence of such elements lying in $D$.

We construct, given a fusion non-Schur element $x \in D$ of the block fusion system $\F_b$, a $p'$-element $z \in C_G(x)$ so that $xz \in G$ is an extra non-Schur element, the projection of whose conjugacy class sum onto $\HH^0(B)$ is non-zero. In particular, this implies the non-vanishing of the BV operator on $\HH^1(B)$, allowing us to prove the following fusion theoretic criterion.

\begin{restatable}{theorem}{MainTheorem}\label{mainthm} Suppose $B=kGb$ is a block of $kG$ with defect group $D$ and with fusion system $\mathcal{F}_b$. If $\mathcal{F}_b$ contains a fusion non-Schur element, then $\operatorname{HH}^1(B)\not=0$.
\end{restatable}

In particular, given a fusion system $\mathcal{F}$ on an abelian $p$-group $D$, each $x \in D \backslash \Phi(D)$ is a fusion non-Schur element, giving Theorem \ref{AbelianTheorem}.

We proceed to show that all fusion systems of blocks at $p>5$ with non-trivial defect admit fusion non-Schur elements. We demonstrate this by first showing a minimal counterexample must be reduction simple, and have a defect group whose conjugates generate the ambient group $G$. A result of Serwene \cite[Lemma 3.4]{Serwene} then implies that should a counterexample exist, one must exist amongst blocks of quasisimple groups with $p'$-center. We proceed to verify this criterion using the classification of finite simple groups.

\begin{restatable}{theorem}{MainTheorem2}\label{mainthm2} Suppose $B$ is a block of $kG$ with non-trivial defect group $D$ and fusion system $\mathcal{F}$ at $p>5$. Then $\mathcal{F}$ contains a fusion non-Schur element.
\end{restatable}

Combining Theorem \ref{mainthm} with Theorem \ref{mainthm2}, we are able to deduce the validity of Conjecture \ref{markus1} in characteristic greater than $5$.

\subsection*{Acknowledgements} The author would like to thank Gabriel Navarro, Markus Linckelmann and Burkhard K{\"u}lshammer for helpful communications. Further, the author would like to thank Radha Kessar and Martin Van Beek for conversations about fusion systems. Finally, the author thanks Ben Briggs and Lleonard Rubio y Degrassi for invaluable conversations and their support.

\subsection*{AI Disclosure} Theorem \ref{mainthm} was proved without the assistance of Generative AI. Theorem \ref{mainthm2} was the result of conversations with ChatGPT Pro; the argument has been verified and rewritten by the author, and the author assumes full responsibility for the correctness of the arguments herein. 

\section{Preliminaries}
\subsection{Blocks}
To begin we collect notation and results already present in the literature. Throughout $G$ will be a finite group and $k$ an algebraically closed field of characteristic $p>0$. $B$ will be a block of $kG$, with defect group $D$ and block idempotent $b$. We often confuse $b$ and $B$ as is standard. We write $EG$ for the subalgebra of $kG$ spanned by central idempotents. All algebras will be algebras over $k$, as will be all unadorned tensor products.

Before beginning, we recall some standard notation. Given an element $z=\sum_{x \in G} a_xx \in kG$, we write $\operatorname{Supp}(z) = \{x \in G: a_x\not=0\}$. Given an element $x \in G$, we can write $x$ uniquely as a product of commuting elements $x=x_px_{p'}$, where $x_p$ (resp. $x_{p'})$ is a $p$-element (resp. $p'$-element). Crucially, $C_G(x)= C_G(x_p)\cap C_G(x_{p'})$. The $p$-section of $x \in G$ consists of those elements $y \in G$ for which $y_p$ is conjugate to $x_p$ in $G$. Given $X \subset G$, we write $\widehat{X} = \sum_{g \in X} g$, and $kX$ for the $k$-span of $X$ in $kG$. Following \cite{AKO}, given $x, g \in G$, we write $x^g := g^{-1}xg$ and $\hspace{0pt}^gx = gxg^{-1}$.

Recall from \cite[Chapter 2]{navarro}, that given a block $B$, we may associate an algebra homomorphism, called a \emph{central character}, which we write $\lambda_b: Z(B) \rightarrow k$ . Abusing notation, we also write $\lambda_b:Z(kG) \rightarrow k$ given by precomposing with the projection from $kG$ to $B$. The next proposition is standard, though crucial to our proof, so we state it thoroughly.

\begin{propositiondefinition}\label{def:defectclass}
Given a block $b$, there exists a conjugacy class $K$ of $p'$-elements of $G$, called a \emph{defect class} of $b$, for which $K \subset \operatorname{Supp}(b)$ and $\lambda_b(\hat{K}  b) \not=0$. Moreover there is $z \in K$ so that $D \in \operatorname{Syl}_p(C_G(z))$.
\end{propositiondefinition}
\begin{proof}
Write $b := \sum_{K \in G/ \sim} a_K \hat{K}$. Then $\lambda_b(b) = \sum_{K \in G/\sim} a_K \lambda_b(\hat{K}b)\not=0$. Therefore, there must be some conjugacy class $K$ of $G$ for which $a_K\not=0$, and $\lambda_b(\hat{K}b) \not=0$. The statement on defect groups is \cite[Corollary 4.5]{navarro}, and that $z$ is a $p'$-element is a consequence of the fact that block idempotents are supported in $p'$-elements, see \cite[(61)]{kulshammer}.
\end{proof}

For $X \subset G$ we write $\pi_X^G:Z(kG)\rightarrow kX$ for the projection map. When $P$ is a $p$-subgroup of $G$, the map $\pi_{C_G(P)}^G$ is an algebra homomorphism called the Brauer map, and written $\Br_P:Z(kG) \rightarrow Z(kC_G(P))$. We need one more result.
\begin{proposition}\cite[(61)]{kulshammer}\label{kulshammerproj} Let $U$ be the $p$-section of $x \in Z(G)$. Then $\pi_U^G$ is $EG$-linear. That is, given $e \in EG$ and $a \in Z(kG)$, $\pi_U^G(ea) = e\pi_U^G(a)$.
\end{proposition}

\subsection{Fusion Systems}

Recall a fusion system $\mathcal{F}$ on $P$ is a category whose objects are subgroups of $P$. Moreover, all morphisms are injective group homomorphisms of the subgroups, each of which factors as an isomorphism in $\mathcal{F}$ followed by an inclusion, containing all conjugation maps by elements of $P$ (see \cite{AKO}).  Given $R, Q \leq P$, we write $Q \sim_{\mathcal{F}} R$ when $Q$ and $R$ are isomorphic in $\mathcal{F}$, and given $x, y\in P$, we write $x \sim_\F y$ when there is an isomorphism in $\mathcal{F}$ sending $x$ to $y$. We will be primarily concerned with \emph{saturated} fusion systems.
\begin{definition}\cite[Part I, Definitions 2.2, 2.4 \& Proposition 2.5]{AKO}

Let $\mathcal{F}$ be a fusion system on a $p$-group $P$.
\begin{itemize}
\item A subgroup $Q \leq P$ is fully-$\mathcal{F}$-centralized (resp. normalized) when for any $R \leq P$ for which $Q \sim_\mathcal{F}R$, we have $|C_P(Q)|\geq|C_P(R)|$ (resp. $|N_P(Q)|\geq|N_P(R)|$).
\item A subgroup $Q \leq P$ is called fully-$\mathcal{F}$-automized if $\operatorname{Aut}_P(Q) \in \operatorname{Syl}_p(\operatorname{Aut}_\mathcal{F}(Q))$.
 \item A subgroup $R \leq P$ is called \emph{receptive} in $\mathcal{F}$ if for each $Q\leq P$, and each $\phi \in \operatorname{Iso}_\mathcal{F}(Q,R)$, writing
$$N_\phi := \{ g \in N_P(Q): \hspace{0pt}^\phi c_g \in \operatorname{Aut}_P(R) \},$$
then there is $\tilde\phi \in \Hom_\F(N_\phi, P)$ so that $\tilde\phi\big|_Q = \phi$. 
\end{itemize}
A fusion system is \emph{saturated} if it has the following properties.
\begin{enumerate}[(i)]
	\item Each subgroup $Q\leq P$ which is fully-$\F$-normalized is also fully-$\mathcal{F}$-centralized and fully-$\F$-automized.
	\item Each subgroup $Q \leq P$ which is fully-$\F$-centralized is also receptive.
\end{enumerate}
\end{definition}
To each finite group $G$ with a chosen Sylow $p$-subgroup, we can associate a fusion system.

\begin{definition}\cite[8.1.2]{linckelmann} Given a group $G$ with Sylow $p$-subgroup $P$, we define $\mathcal{F}_P(G)$ to be the fusion system on $P$ with morphisms $Q \xrightarrow\varphi R$, for which there exists $g\in G$ so that $\varphi(x) = x^g$ for all $x\in Q$.
\end{definition}

We call a fusion system $\mathcal{F}$ \emph{realisable} when it is isomorphic to the fusion system of a group.

A Brauer pair on $kG$ is a pair $(Q,f)$ where $Q$ is a $p$-subgroup of $G$ and $f$ is a block of $kC_G(Q)$. A $(G,B)$-Brauer pair is a Brauer pair $(Q,f)$ so that $\operatorname{Br}_Q(b)f =f$. There is a natural action of $G$ on $(G,B)$-Brauer pairs, which we write $(Q,f) \mapsto \hspace{0pt}^x(Q,f)$, and a notion of inclusion of Brauer pairs on $kG$, written $(Q,f) \leq (P,e)$; for details see \cite[6.3]{linckelmann}. A $(G,B)$-Brauer pair maximal under inclusion is called a maximal $(G,B)$-Brauer pair. We can choose $(D,e_D)$ a maximal $(G,B)$-Brauer pair for some idempotent $e_D$ of $kC_G(D)$.
\begin{proposition}\cite[6.3.3, 6.3.6]{linckelmann}\label{uniqueidemp}
Let $(P,e_P)$ be a Brauer pair on $kG$, and take a subgroup $Q$ of $P$. There is a unique block $B_Q$ of $kC_G(Q)$, with idempotent $e_Q$ so that $(Q,e_Q) \leq (P,e_P)$. Moreover, $(P,e_P)$ is a $(G,B)$-Brauer pair if and only if $(Q,e_Q)$ is.
\end{proposition}
We can now define the fusion system associated to a block.
\begin{definition}\cite[8.5.1]{linckelmann}
    Let $(D,e_D)$ be a maximal $(G,B)$-Brauer pair. For each $Q \leq D$, write $(Q,e_Q)$ for the unique Brauer pair contained in $(D,e_D)$. We define $\mathcal{F}_{(D,e_D)}(G,B)$ to be the category with objects the subgroups of $D$, with morphisms $\varphi:Q \rightarrow R$ for which there exists $g \in G$ so that $\varphi(x)= \hspace{0pt}^gx$ for all $x \in Q$, and $\hspace{0pt}^g(Q,e_Q) \leq (R, e_R)$.
\end{definition}

The fusion systems of groups and of blocks under consideration are saturated, as we have assumed that $k$ is algebraically closed. 

 \begin{proposition}\cite[8.1.4, 8.5.3]{linckelmann}\label{goodcentral} Let $P$ denote a Sylow $p$-subgroup of $G$. We can describe fully centralized subgroups of $\mathcal{F}_P(G)$ and $\mathcal{F}_{(D,e_D)}(G,B)$ as follows.
     \begin{enumerate}[(i)]
    \item $Q \leq P$ is fully-$\mathcal{F}_P(G)$-centralized if and only if $C_P(Q)$ is a Sylow $p$-subgroup of $C_G(Q)$.
    \item $Q \leq D$ is fully-$\mathcal{F}_{(D,e_D)}(G,B)$-centralized if and only if $C_D(Q)$ is a defect group of $kC_G(Q)e_Q$. 
     \end{enumerate}
 \end{proposition}  
 \begin{definition}\cite[4.26]{Craven} Given a fusion system $\mathcal{F}$ on $P$, and a subgroup $Q \leq P$, the category $C_\mathcal{F}(Q)$ has objects the subgroups of $C_{P}(Q)$, and has morphisms $\phi \in \Hom_{C_\F(Q)}(R,S)$ when $\phi \in \Hom_\F(R,S)$ extends to $\bar{\phi} \in \Hom_\F(QR, QS)$ so that $\bar{\phi}|_{Q} = \operatorname{Id}_Q$. 
 \end{definition}
 \begin{proposition}\cite[8.5.4, 8.3.11(ii)]{linckelmann}\label{prop:centralizerblock} Let $P$ denote a Sylow $p$-subgroup of $G$, and $(D,e)$ be a maximal $(G,B)$-Brauer pair for $B$ a block of $kG$.
 \begin{enumerate}[(i)]
 \item  Suppose $Q$ is a fully-$\mathcal{F}_{P}(G)$-centralized subgroup of $P$. Then $C_{\mathcal{F}_P(G)}(Q) = \mathcal{F}_{C_P(Q)}(C_G(Q))$.
 \item  Suppose $Q$ is a fully-$\mathcal{F}_{(D,e)}(G,B)$-centralized subgroup of $D$. Then $(C_D(Q), e_{QC_D(Q)})$ is a maximal $(C_G(Q),e_Q)$-Brauer pair and we have $C_{\mathcal{F}_{(D,e)(G,B)}}(Q) = \mathcal{F}_{(C_D(Q), e_{QC_D(Q)})}(C_G(Q),e_Q)$.  
\end{enumerate}
 \end{proposition}
 Our eventual aim will be to reduce the non-vanishing question to a property of the blocks of quasisimple groups. As such, we need an analogue of normal subgroups for fusion systems.
  \begin{definition}\cite[Definition 4.55]{Craven} Take $Q$ a subgroup of $P$. Then $Q$ is called:
 \begin{itemize}
     \item Weakly $\mathcal{F}$-closed if whenever $\phi \in \Hom_\F(Q,P)$, $\phi(Q)=Q$. 
     \item Strongly $\mathcal{F}$-closed if whenever $R$ is a subgroup of $Q$, and $\phi \in \Hom_\F(R,P)$, then $\phi(R) \leq Q$. 
     \end{itemize}
      A fusion system is called reduction simple if it admits no non-trivial proper strongly closed subgroups.
  \end{definition}
 \begin{definition}\cite[Definition 5.9]{Craven}
    Let $Q \trianglelefteq P$, and write $q: P \rightarrow P/Q$ for the quotient map. The factor system $\mathcal{F}/Q$ is the category with objects subgroups of $P/Q$, and for any two subgroups $R$ and $S$ containing $Q$, we define $\Hom_{\F/Q}(R/Q, S/Q)$ as the set of homomorphisms $\tilde\phi$, so that there exists $\phi \in \Hom_\F( R, S)$, with $\phi(Q) = Q$ and $q \circ \phi = \tilde \phi \circ q$. 
 \end{definition}
 In fact for $\mathcal{F}$ saturated, when $Q$ is weakly $\mathcal{F}$-closed $\mathcal{F}/Q$ is saturated, by \cite[Proposition 5.10]{Craven}.
 \begin{definition}\cite[Definition 5.1]{Craven}
    Given two fusion systems $\mathcal{E}, \mathcal{F}$ on $p$-groups $P$ and $Q$ respectively, a morphism of fusion systems $ \mathcal{E} \rightarrow \mathcal{F}$ consists of a group homomorphism $\alpha: P \rightarrow Q$ and a functor $\Phi:\mathcal{E} \rightarrow \mathcal{F}$ obeying the following properties.
    \begin{itemize}
        \item For each $R \leq P$, $\alpha(R) = \Phi(R)$. 
        \item For all  $R, T \leq P$, and each morphism $\phi \in \Hom_{\mathcal{E}}(R,T)$, $\Phi(\phi) \circ \alpha = \alpha \circ \phi$.
    \end{itemize}
    We call $\alpha$ the underlying group homomorphism of $\mathcal{E} \rightarrow \mathcal{F}$.
 \end{definition}
 Note that in the above, the functor $\Phi$ is entirely determined by $\alpha$.

 \section{Fusion non-Schur elements and non-vanishing}
In this section, we will prove our fusion theoretic criterion for non-vanishing.  Given a $p$-group $P$ we write $\Phi(P)$ for its Frattini subgroup. For a finite group $G$ we write $E^p(G) = G^p[G,G]$, where $G^p$ is the subgroup generated by $p^{th}$-powers of elements of $G$. Note that $E^p(P) = \Phi(P)$, and for $x \in G$, we have $x \in E^p(G)$ if and only if $x \in \operatorname{ker}(\varphi)$ for each group homomorphism $\varphi:G \rightarrow k^+$. We restate the definition of fusion non-Schur elements for convenience.

\MainDef*

In realisable fusion systems on non-trivial $p$-groups at $p>5$, these elements always exist. This bound is sharp, as for realisable fusion systems, they agree up to conjugacy with extra strong non-Schur elements of \cite{BR}, where explicit counterexamples exist for $p\leq 5$. 
 \begin{theorem}\label{thm:realisable}
 Suppose $p> 5$, and $G$ is a finite group with non-trivial Sylow $p$-subgroup $P$. Then $\mathcal{F}_P(G)$ admits a fusion non-Schur element.
 \end{theorem}
 \begin{proof} 
	 From \cite[Theorem 5.1]{BR} we know that for $p>5$, each group $G$ as above has an extra-strong non-Schur element. This is a $p$-element $y \in G$ which has $y \not \in E^p(C_G(y))$. We may choose a Sylow $p$-subgroup $Q$ so that $C_Q(y)$ is a $p$-Sylow of $C_G(y)$; then $Q$ is conjugate to $P$ by some $h\in G$. Writing $x := y^h$, we have $x \in P$, and know that $C_P(x)$ is a Sylow $p$-subgroup of $C_G(x)$. Thus by Proposition \ref{goodcentral}, we know that $\langle x \rangle \leq P$ is fully $\mathcal{F}_P(G)$-centralized. Utilising \cite[Lemma 3.8]{BR} we have $x \in \Phi(C_P(x))$ if and only if $x \in E^p(C_G(x))$. As the property of being extra non-Schur is invariant under conjugation, and as $x \not \in E^p(C_G(x))$, we see $x$ is a fusion non-Schur element.
 \end{proof}
 In the same vein as \cite{BR}, our proof relies on non-vanishing of the Batalin-Vilkovisky (BV) operator on $\operatorname{HH}^1(B)$, defined in \cite{tradler}.

\begin{lemma}\label{bvnonvan} Suppose $x \in G$ has $x \not \in E^p(C_G(x))$, and write $K := x^G$. If $b\hat{K} \not=0 $, then the BV operator $\Delta_B^{1} \colon\operatorname{HH}^1(B) \rightarrow \HH^0(B)$, is non-zero. In particular $\HH^1(B)\not=0$.
\end{lemma}
\begin{proof}
    As $x \not \in E^p(C_G(x))$ there is a group homomorphism $\varphi: C_G(x) \rightarrow k^+$ so that $\varphi(x) \not=0$. Following \cite[Proposition 2.3]{BR} we write $\gamma_x :\operatorname{H}^1(C_G(x),k) \rightarrow \HH^1(kG)$ for the inclusion induced by the centralizer decomposition of Hochschild cohomology. Then from \cite[Theorem 1.1 (ix)]{BKL} and \cite[Proposition 2.3]{BR} we know the Batalin-Vilkovisky operator $\Delta_{kG}^1:\operatorname{HH}^1(kG) \rightarrow \operatorname{HH}^{0}(kG)$ is given by $\Delta_{kG}^1(\gamma_x(\varphi)) = \varphi(x)\hat{K}$. Further from \cite[Proposition 2.6]{BR}, we know that with natural choice of form on $B$, we have: $$\Delta_B^1( \gamma_x (\varphi)b) = \Delta_{kG}^1(\gamma_x(\varphi))b = \varphi(x)\hat{K}b \not=0.$$
    Thus, $\Delta_B^1(\gamma_x(\varphi)b) \not=0$, finishing the proof.
\end{proof}
We need one more proposition before we prove our theorem. 
\begin{proposition}\label{phipep} Suppose $x \in G$ is central, and $Q \in \operatorname{Syl}_p(G)$. Then $x_p \in \Phi(Q) \iff  x\in E^p(G)$.
\end{proposition}
\begin{proof}
   We see that $x_p$ is central in $G$ also, so by \cite[Lemma 3.8]{BR} we know that $x_p \in E^p(Q)$ if and only if $x_p \in E^p(G)$. But given $\varphi:G \rightarrow k^+$, we have $\varphi(x_p) = \varphi(x)$, so that $x \in E^p(G)$ if and only if $x_p \in E^p(G)$. As $E^p(Q) = \Phi(Q)$, we are done.
\end{proof}

\MainTheorem*
 \begin{proof}[Proof of Theorem \ref{mainthm}]
Suppose $x \in D$ is a fusion non-Schur element of $\mathcal{F}_{(D,e_D)}(G,B)$. Then by Proposition \ref{uniqueidemp} we have a block $B_x$ of $kC_G(x)$ with block idempotent $e_x$, so that $(\langle x \rangle, e_x) \leq (D,e_D)$. Now, as $\langle x \rangle$ is fully centralized, we know that $B_x$ has defect $C_D(x)$ by Proposition \ref{goodcentral}. Take $K$ a defect class of $B_x$ in $C_G(x)$, and choose $z \in K$ so that $C_D(x) \in \operatorname{Syl}_p(C_{C_G(x)}(z))$. We claim that $xz \not \in E^p(C_G(xz))$ and writing $L:= (xz)^G$, that $b\hat{L} \not=0$. Having proved these claims, we are done by Lemma \ref{bvnonvan}.

We first show that $xz \not \in E^p(C_G(xz))$. We note that $x,z$ commute as $z \in C_G(x)$, and $x$ (resp $z$) is a $p$-element (resp $p'$-element). So $(xz)_p = x$ and $(xz)_{p'} = z$. Thus, $C_G(xz) = C_G(x) \cap C_G(z) = C_{C_G(x)}(z)$. Now by our choice of $z$, $C_D(x)$ is a Sylow $p$-subgroup of $C_{C_G(x)}(z) = C_G(xz)$. Then using Proposition \ref{phipep} we know that $ x \in \Phi(C_D(x))$ if and only if $ xz \in E^p(C_G(xz))$. Then, as $x$ is a fusion non-Schur element, we are done.

Write $U$ for the $p$-section of $x$ in $C_G(x)$. We now demonstrate that $b\hat{L} \not =0$, by showing that the element $\alpha:= \lambda_{e_x}(e_x \pi_U^{C_G(x)}(\Br_{\langle x \rangle}(b \hat{L})))$ is non-zero.

We first note that as $\operatorname{Br}_{\langle x \rangle}$ is an algebra homomorphism, $\operatorname{Br}_{\langle x \rangle}(b\hat{L}) = \Br_{\langle x \rangle}(b)\operatorname{Br}_{\langle x \rangle}(\hat{L})$, and $\Br_{\langle x \rangle}(b)$ is an idempotent. Moreover from Proposition \ref{kulshammerproj} we know that $\pi_U^{C_G(x)}$ is $EC_G(x)$-linear. Hence we see that $\pi_U^{C_G(x)}(\Br_{\langle x \rangle}(b \hat{L})) = \Br_{\langle x \rangle}(b)\cdot \pi_U^{C_G(x)}(\Br_{\langle x \rangle}(\hat{L}))$. For convenience, write $\beta := \pi_U^{C_G(x)}(\Br_{\langle x \rangle}(\hat{L}))$. 

Finally as $(\langle x \rangle, e_{x})$ is a $(G,B)$-Brauer pair by Proposition \ref{uniqueidemp}, we have that $\operatorname{Br}_{\langle x \rangle}(b)e_{x} = e_{x}$. In total $\alpha = \lambda_{e_x}(\beta e_{x})$. So we must understand the element $\beta \in Z(kC_G(x))$. By definition, we have: 
$$\operatorname{Br}_{\langle x \rangle}(\hat{L})\; = \hspace{-10pt}\sum_{\substack{g\,\in\,\rightcosets{C_G(xz)}{G}\\ (xz)^g\,\in\,C_G(x)}}  \hspace{-10pt}x^gz^g.$$

We seek to apply $\pi_U^{C_G(x)}$ to this sum, where $U$ is the $p$-section of $x$ in $C_G(x)$. Notably, $x$ is central in $C_G(x)$ so that $U = \{ y \in C_G(x): y_p = x \}$. Moreover, $(x^gz^g)_p=x^g$, so that any summand not eliminated by $\pi_U^{C_G(x)}$ must have $x^g=x$, that is $g \in C_G(x)$. Therefore, recalling $K = z^{C_G(x)}$, we have:
$$\beta = \hspace{-10pt}\sum_{g\,\in\,\rightcosets{C_G(xz)}{C_G(x)}} \hspace{-15pt} xz^g = x\hat{K}$$
Thus we have $\alpha = \lambda_{e_x}(x\hat{K} e_{x}) = \lambda_{e_x}(x e_{x})\lambda_{e_x}(\hat{K}e_{x})$. This is non-zero as $xe_{x}$ is a unit in $Z(B_x)$, and $\lambda_{e_x}(\hat{K}e_{x})\not=0$ by definition of defect classes. 
 \end{proof}
We immediately obtain a result valid for all primes.
\begin{corollary}\label{cor:abelian} Let $B$ be a block of $kG$ with defect $D$. Each $x \in Z(D) \cap (D \backslash \Phi(D))$ is a fusion non-Schur element of $\mathcal{F}_b$. As such, if $Z(D) \not \subset \Phi(D)$ then  $\HH^1(B)\not=0$. \end{corollary}
\begin{proof} Take any $x \in Z(D)$ for which $x \not \in \Phi(D)$. Then $x$ is fully centralized, and $x \not \in \Phi(C_D(x)) = \Phi(D)$, so $x$ is a fusion non-Schur element.
\end{proof}
In particular, we obtain Theorem \ref{AbelianTheorem}.
\AbelianTheorem*
\section{Reduction to blocks of Quasisimple groups.}
 \subsection{Reduction results for saturated fusion systems.}
 We now seek to prove the existence of fusion non-Schur elements in block fusion systems at $p>5$. We achieve this by proving a reduction theorem to blocks of quasisimple groups.  Throughout $\mathcal{F}$ will be a fusion system on a $p$-group $P$. We first prove some general reduction results for saturated fusion systems. 
 \begin{lemma}\label{lem:subsystem} Let $\mathcal{E}$ be a fusion subsystem of $\mathcal{F}$ on the same $p$-group $P$. Then $\mathrm{NS}(\mathcal{F}) \implies \mathrm{NS}(\mathcal{E})$.
 \end{lemma}
 \begin{proof} Take $x \in P$ a fusion non-Schur element of $\mathcal{F}$. We know that $\langle x \rangle$ is fully-$\mathcal{F}$-centralized, and that $ x \not \in \Phi(C_P(x))$. As $\mathcal{E}$ is a fusion subsystem of $\mathcal{F}$, $\langle x \rangle$ is surely fully-$\mathcal{E}$-centralized, so that $x$ is a fusion non-Schur element of $\mathcal{E}$.
 \end{proof}
 \begin{lemma}\label{lem:surjection} Suppose $\mathcal{F}$, $\mathcal{E}$ are saturated fusion systems on $P,Q$ respectively, and $\Phi: \mathcal{F} \rightarrow \mathcal{E}$ is a morphism of fusion systems, with underlying homomorphism $\alpha: P \rightarrow Q$.  If $\alpha$ is surjective then $\mathrm{NS}(\mathcal{E}) \implies \mathrm{NS}(\mathcal{F})$.
 \end{lemma}
 \begin{proof} 
     Take $x \in P$ so that $\alpha(x) \in Q$ is a fusion non-Schur element of $\mathcal{E}$. Thus $\langle \alpha(x) \rangle$ is fully-$\mathcal{E}$-centralized and there exists a group homomorphism $\psi:C_Q(\alpha(x)) \rightarrow k^+$ so that $\psi(\alpha(x)) \not=0$. Take $y \in P$ so that $y \sim_\mathcal{F} x$, by $\theta \in \operatorname{Iso}_\F(\langle y \rangle, \langle x \rangle)$, with $\langle y \rangle$ fully-$\mathcal{F}$-centralized. Then $\alpha(y) \sim_\mathcal{E} \alpha(x)$ by the morphism $\Phi(\theta) \in \operatorname{Iso}_\mathcal{E}( \langle \alpha(y) \rangle, \langle \alpha(x) \rangle)$. By saturation of $\mathcal{E}$, $ \langle \alpha(x) \rangle$ is receptive, and as such there exists a morphism $C_Q(\alpha(y)) \rightarrow C_Q(\alpha(x))$ sending $\alpha(y)$ to $\alpha(x)$. Now the composition: 
      $$C_P(y) \xrightarrow{\alpha} C_{Q}( \alpha(y))\rightarrow C_{Q}(\alpha(x)) \xrightarrow\psi k^+$$
      demonstrates that $y$ is a fusion non-Schur element of $\mathcal{F}$.
 \end{proof}
 \begin{proposition}\label{prop:quotient} Suppose $\mathcal{F}$ is a saturated fusion system and $Q$ is a strongly $\mathcal{F}$-closed subgroup. Then $\mathrm{NS}(\mathcal{F}/Q) \implies \mathrm{NS}(\mathcal{F})$.
 \end{proposition}
 \begin{proof} By \cite[Proposition 5.7, Theorem 5.14]{Craven}, there exists a morphism of fusion systems of the form $\mathcal{F} \rightarrow \mathcal{F}/Q$ meeting the criteria of Lemma \ref{lem:surjection}.
 \end{proof}
 \begin{definition}
	 Given fusion systems $\mathcal{E}, \mathcal{F}$ on the same $p$-group $P$, we say that $\mathcal{F}$ is \emph{wrapped} by $\mathcal{E}$ when the following factorisation condition holds: for each $Q \leq P$ and $\phi \in \Hom_\F(Q,P)$, there are $\psi \in \Hom_{\mathcal{E}}(Q,P)$ and $\alpha \in \operatorname{Aut}(P)$ so that $\phi = \alpha \circ \psi$.
\end{definition}
    The next two results will be crucial in our reduction to blocks of quasisimple groups.
 \begin{lemma}\label{lem:automorph}
	 Let $\mathcal{E}, \mathcal{F}$ be fusion systems on the same group $P$, with $\mathcal{F}$ wrapped by $\mathcal{E}$. Then we have that $\mathrm{NS}(\mathcal{E}) \implies \mathrm{NS}(\mathcal{F})$.
 \end{lemma}
 \begin{proof}
     Take $x \in P$ a fusion non-Schur element of $\mathcal{E}$. Then take $y \in P$ and $\phi \in \Hom_\F(\langle x \rangle, P)$ with $\phi(x) =y$, so that $\langle y \rangle$ is fully-$\mathcal{F}$-centralized. We may write $\phi = \alpha \circ \psi$ as above, wherein $|C_P(x)| \geq |C_P(\psi(x))|$ as $x$ is fully-$\mathcal{E}$-centralized, and $|C_P(\psi(x))| = |C_P(y)|$ as $\operatorname{Aut}(P)$ preserves the order of centralizers. Thus $|C_P(x)| = |C_P(y)|$, so $x$ is fully-$\mathcal{F}$-centralized and $x$ is a fusion non-Schur element of $\mathcal{F}$.
 \end{proof}
 
 \begin{lemma}\label{lem:stronglyclosed} Suppose $p > 5$ and $\mathcal{F}$ is a saturated fusion system for which $\mathrm{NS}(\mathcal{F})$ is false. If $Q$ is a strongly $\mathcal{F}$-closed subgroup then $P= QC_P(Q)$. 
 \end{lemma}
 \begin{proof}
    We first note that as $Q$ is strongly $\mathcal{F}$-closed and so fully-$\mathcal{F}$-normalized, by saturation, $Q$ is fully-$\F$-automized; that is, we have $\operatorname{Aut}_P(Q) \in \operatorname{Syl}_p(\operatorname{Aut}_\mathcal{F}(Q))$. As such, we know that $\operatorname{Out}_P(Q)$  is a Sylow $p$-subgroup of $\operatorname{Out}_\mathcal{F}(Q)$, and that $\mathcal{F}_{\operatorname{Out}_P(Q)}(\operatorname{Out}_\mathcal{F}(Q))$ is a saturated fusion system.
    
    As $Q$ is strongly $\mathcal{F}$-closed, it is certainly normal in $P$, and so writing $c_u$ for the conjugation action of $u\in P$ on $Q$, the morphism $P/Q \xrightarrow{\alpha} \operatorname{Out}_P(Q)$ sending $uQ \mapsto c_u\operatorname{Inn}(Q)$ is a well defined group homomorphism. We seek to show that this gives rise to a morphism $\Phi:\mathcal{F}/Q \rightarrow\mathcal{F}_{\operatorname{Out}_P(Q)}(\operatorname{Out}_\mathcal{F}(Q))$. To do this, we show that given $Q \leq R, T \leq P$, and a homomorphism $\tilde\psi \in \operatorname{Hom}_{\F/Q}(R/Q,T/Q)$ that the morphism $\Phi(\tilde\psi)$ defined by $\Phi(\tilde\psi) \circ \alpha = \alpha \circ \tilde\psi$ is equal to conjugation by an element of $\operatorname{Out}_\mathcal{F}(Q)$.

    Taking subgroups, $Q \leq R, T \leq P$ and a morphism $\tilde\psi$ as above, we have $\psi \in \operatorname{Hom}_{\mathcal{F}}(R,T)$ giving rise to $\tilde\psi$ by the definition of the factor fusion system. Taking $u \in R$ we have $\Phi(\tilde\psi)(c_{u}\operatorname{Inn}(Q)) = c_{\psi(u)}\operatorname{Inn}(Q)$. But notably taking $x \in Q$, we have:
    $$ c_{\psi(u)}(x) = \psi(u)x \psi(u)^{-1} =  \hspace{0pt}^{\psi|_{Q}}(c_u)(x).$$
    That is, $\Phi(\tilde\psi)$ is given by conjugating $c_u \in \operatorname{Out}_P(Q)$ by $\psi|_{Q}\operatorname{Inn}(Q) \in \operatorname{Out}_{\mathcal{F}}(Q)$. Supposing that $\operatorname{Out}_P(Q)$ is non-trivial, by Lemma \ref{lem:surjection} and Theorem \ref{thm:realisable} we see that $\mathcal{F}/Q$ must have a fusion non-Schur element. But by Proposition \ref{prop:quotient} this implies $\mathcal{F}$ has a fusion non-Schur element, so then $\operatorname{Out}_P(Q)$ is trivial. This tells us that $P= C_P(Q)Q$.
 \end{proof}
 \subsection{Reduction for fusion systems of blocks.}
 We now prove the reduction to non-trivial fusion systems of blocks of quasisimple groups with $p'$-center. Throughout $G$ is a finite group, and $b$ is a block of $kG$ of defect $D$. We write $(D,e)$ for a maximal $(G,b)$-Brauer pair and $\mathcal{F}_b := \mathcal{F}_{(D,e)}(G,b)$ for the fusion system of the block. We first state some results in the literature. 
 \begin{proposition}\label{prop:centralquotient} Take $Z \leq G$ a central $p$-subgroup, write $\bar{G} = G/Z$, and $\bar{b}$ for the block of $k\bar{G}$ dominated by $b$. Then we have that $Z \leq D$, that $D/Z$ is the defect group of $\bar{b}$, and that there exists a maximal $(\bar{G}, \bar{b})$-Brauer pair $(D/Z,f)$ such that $\mathcal{F}_b/Z\cong \mathcal{F}_{(D/Z, f)}(\bar{G}, \bar{b})$.
 \end{proposition}
 \begin{proof} From \cite[Proposition 4.2]{KLR} we know that $G \rightarrow \bar{G}$ induces a surjective on objects and full morphism of fusion systems $\mathcal{F}_{(D,e)}(G,b)/ Z \rightarrow \mathcal{F}_{(D/Z, f)}(\bar{G}, \bar{b})$. Then by the first isomorphism theorem for fusion systems (see \cite[Corollary 5.15]{Craven}) we are done.
 \end{proof}
 \begin{lemma}\label{prop:serwene}\cite[Lemma 3.4]{Serwene}
 Take $\mathcal{F}$ the fusion system of a block of group $G$ with non-abelian defect group $D$. If $G  = \langle D^g: g \in G \rangle$ and $\mathcal{F}$ is reduction simple, then there is a quasisimple group $L$ with $p'$-center, admitting a block with fusion system $\mathcal{F}$.
 \end{lemma}
 
 \begin{proposition}\label{prop:usefulprop}
	 Take $N \leq G$ a normal subgroup, and $b$, $\mathcal{F}_b$ as above. 
    \begin{enumerate}[(i)]
        \item\cite[Lemma 6.1]{KS}\label{part:usefulprop1} There exists a block $c$ of $kN$ so that $c$ is $D$-stable with $cb \not=0$ and with defect group $D \cap N$.
        \item\cite[Part IV, Proposition 6.3]{AKO}\label{part:usefulprop2} Take $d$ a block of $kN$, and suppose $b$ covers $d$. Write $I$ for the stabilizer of $d$ in $G$. Then there is a block $e$ of $kI$ covering $d$ so that $e$ and $b$ have isomorphic fusion systems.
	\item\cite[Theorem 3.5]{KS}\label{part:normalsubsystem} Suppose $D \leq N$, take $d$ a $G$-stable block of $kN$ such that $db\not=0$, and write $\mathcal{F}_d$ for the fusion system of the block $d$. Then there exists a fusion system $\mathcal{H}$ on $S \leq G$ so that $\mathcal{F}_b$ and $\mathcal{F}_d$ are both fusion subsystems of $\mathcal{H}$ on $D$, with $\mathcal{F}_d$ a weakly normal subsystem. In particular, $\mathcal{F}_b$ is wrapped by $\mathcal{F}_d$.
	\item\cite[Lemma 3.2]{Eisele}\label{part:quotients} Suppose $N$ is a central $p'$-subgroup of $G$, and take $\bar{b}$ a block of $kG/N$ of defect $\bar{D}$. Then there exists a unique block $b$ of $kG$ of defect $D$, whose quotient maps onto $\bar{b}$. Further, $D \cong \bar{D}$, and $\mathcal{F}_{\bar{b}}$ is a subsystem of $\mathcal{F}_{b}$ on the same $p$-group.
	\item\label{part:normalimplication} If $D \leq N$, and $c$ is a block of $kN$ covered by $b$, then $\mathcal{F}_b$ is wrapped by $\mathcal{F}_c$, so that in fact $\mathrm{NS}(\mathcal{F}_c) 
\implies \mathrm{NS}(\mathcal{F}_b)$. 
    \end{enumerate}
 \end{proposition}
 \begin{proof} $(i)$-$(ii)$ are proved in their respective references. For part $(iii)$, note that weak normality gives the Frattini condition (see \cite[Part I, Definition 6.1]{AKO}), which tells us the following. For each $R \leq D$ and $\psi \in \Hom_{\mathcal{F}_b}(R, D)$ we have $\alpha \in \operatorname{Aut}(D)$ and $\psi_0 \in \Hom_{\mathcal{F}_d}(R,D)$ so that $\psi = \alpha \circ \psi_0$. As such $\mathcal{F}_b$ is wrapped by $\mathcal{F}_d$. Part $(iv)$ is seen by specialising the result of \cite{Eisele} when $N$ is a $p'$-central subgroup, and using \cite[Theorem 5.8.8]{Nagao} to verify their criteria. It remains to prove part $(v)$. By part $(ii)$, we may assume without restriction of generality that $c$ is $G$-stable. By part $(iii)$, we know that $\mathcal{F}_b$ is wrapped by $\mathcal{F}_c$. Then Lemma \ref{lem:automorph} completes the claim.  
\end{proof}
\begin{proposition}\label{prop:reductsimple} Assume $p>5$, and that there exists a block fusion system not admitting a fusion non-Schur element. Let $\mathcal{F}_b$ be such a fusion system with $|D|$ minimal amongst non-trivial examples. Then $\mathcal{F}_b$ is reduction simple.
\end{proposition}
\begin{proof}
    Suppose $Q \leq D$ is a non-trivial proper strongly closed subgroup. Then by Lemma \ref{lem:stronglyclosed} we know that $D=QC_D(Q)$. We will deduce the existence of a fusion non-Schur element by constructing two fusion systems on $D/Q$, and then using Lemma \ref{lem:automorph} and Proposition \ref{prop:quotient}. 
    
    We first note that all $\mathcal{F}_b$-automorphisms of $Q$ extend to $\mathcal{F}_b$-automorphisms of $D$. As $Q$ is strongly $\mathcal{F}_b$-closed it is certainly fully centralized. As such, saturation of $\mathcal{F}_b$ requires that for any $\phi \in \operatorname{Aut}_{\mathcal{F}_b}(Q)$, $\phi$ extends to a homomorphism $\bar\phi:N_\phi \rightarrow D$, where:
    $$N_\phi := \{g \in N_D(Q): \hspace{0pt}^\phi c_g \in \operatorname{Aut}_D(Q)\}.$$
    We note that $N_D(Q) =D$, and that for all $g \in D$, there is $h \in Q$ with $c_g = c_h$. Then as $\operatorname{Inn}(Q)$ is a normal subgroup of $\operatorname{Aut}_{\F_b}(Q)$, we have $N_\phi = D$.
    
   Recall that as $Q$ is fully-$\mathcal{F}_b$-centralized, by Proposition \ref{prop:centralizerblock} we know that $\mathcal{E}:=C_{\F_b}(Q)$ is the fusion system of a block $b_Q$ of $C_G(Q)$. Now, $Z: =Z(Q)$ is a central subgroup of $C_G(Q)$, so by Proposition \ref{prop:centralquotient}, $\mathcal{E}/Z$ is the fusion system of a block of $C_G(Q)/Z$, with defect group $C_D(Q)/Z$. As such by minimality of $|D|$, we have $\mathrm{NS}(\mathcal{E}/Z)$.

   Writing $R := C_D(Q)$ we note $R/Z \cong D/Q$; as such $\mathcal{F}_b/Q$ and $\mathcal{E}/Z$ are fusion systems on the same group. Indeed, take $Q \leq U \leq D$, and $\tilde\psi \in \operatorname{Hom}_{\mathcal{F}_b/Q}(U/Q, D/Q)$, induced by $\psi: U \rightarrow D$ with $\psi(Q) = Q$. By the above, $\psi{\big|_Q}$ extends to $\beta \in \operatorname{Aut}_{\F_b}(D)$. Writing $\theta = \beta^{-1}\psi$, we see that $\theta$ fixes $Q$ pointwise. Moreover, as $Q \leq U$, we see that $Q(U \cap R) = U$, meaning that $\theta|_{U \cap R}$ extends to the map $\theta \in \Hom_{\F_b}(Q(U \cap R), D)$, and thus $\theta\big|_{U \cap R} \in \operatorname{Hom}_\mathcal{E}(U\cap R, R)$. Then \cite[Proposition 5.16]{Craven} ensures that the image of $\mathcal{E}$ in $\F_b/Q$ is the same as the image of $\mathcal{E}/Z$, so we have $\tilde\theta \in \operatorname{Hom}_{\mathcal{E}/Z}(U / Q, D/Q)$ and $\tilde\beta \in \operatorname{Aut}(D/Q)$ so that $\tilde\psi = \tilde\beta \circ \tilde\theta$. As such, $\mathcal{F}_b/Q$ is wrapped by $\mathcal{E}/Z$. Thus $\mathrm{NS}(\mathcal{E}/Z)$ implies $\mathrm{NS}(\mathcal{F}_b/Q)$ by Lemma \ref{lem:automorph}, which in turn implies $\mathrm{NS}(\F_b)$ by Proposition \ref{prop:quotient}.
\end{proof}
\begin{proposition}\label{prop:quasisimple} Take $p>5$. Suppose that there exists a $p$-block of non-trivial defect whose fusion system $\mathcal{F}_b$ admits no fusion non-Schur elements. Then there is a $p$-block with non-trivial defect of a quasisimple group with $p'$-center whose fusion system admits no fusion non-Schur elements.
\end{proposition}
\begin{proof}
    Choose $G$ a group with a block $b$ as hypothesized, first so that $|D|$ is minimal, and second with $|G|$ minimal. We note that $D$ is not abelian by Corollary \ref{cor:abelian}, and by Proposition \ref{prop:reductsimple}, $\mathcal{F}_b$ is reduction simple. 
     Write $N:= \langle D^g: g \in G\rangle$.  Suppose that $N$ is a proper subgroup of $G$. Then by Proposition \ref{prop:usefulprop} (\ref{part:usefulprop1}) there exists a block $c$ of $kN$ covered by $b$, which by minimality of $G$ has $\mathrm{NS}(\mathcal{F}_c)$. But by Proposition \ref{prop:usefulprop} (\ref{part:normalimplication}), as $D \leq N$, this implies $\mathrm{NS}(\mathcal{F}_b)$, a contradiction. So we must have $G = N$, and applying Lemma \ref{prop:serwene} we are done.  
\end{proof}
 \section{Existence of fusion non-Schur elements for sporadic and alternating groups.}
 In this section we verify that the fusion systems of blocks of quasisimple groups admit fusion non-Schur elements, when the underlying simple group is sporadic or alternating. Recall that a group $E$ is called quasisimple when it is a perfect central extension of a simple group; that is $E = [E,E]$ and $E/Z(E)$ is simple. Recall that to every perfect group $G$, there exists a universal perfect central extension $\hat{G}$; it is universal in the sense that every other perfect central extension arises as a factorisation of $\hat{G} \rightarrow G$. The kernel of $\hat{G} \rightarrow G$ is called the Schur Multiplier of $G$, which we write $M(G)$. Throughout, $b$ is a block of $kG$ of non-trivial defect $D$, with fusion system $\mathcal{F}_b$. Given a set $X \subset \mathbb{N}$, and a collection of primes $\mathcal{P}$, we write $|X|_{\mathcal{P}}$  for the greatest divisor of $|X|$ whose prime factors lie in $\mathcal{P}$. In particular, when $\mathcal{P} = \{p\}$, we write $|X|_p$ for the $p$-part of $|X|$.

 Before the next result, we inductively define a family of $p$-groups that will arise frequently as defect groups. Set $W_{0} = C_{p}$, and $W_{t} = W_{t-1} \wr C_p$.
 \begin{proposition}\label{prop:guarantee} The following criteria guarantee that $\mathcal{F}_b$ contains a fusion non-Schur element.
	 \begin{enumerate}[(i)]
	 \item\label{part:maximal} We have $p>5$ and $b$ is a block of maximal defect in $kG$.
	 \item\label{part:lessthanp3} We have $p>5$ and $b$ is a block of $kG$ where $|G|_p \leq p^3$. 
	 \item\label{part:lowerexponent} We have $\operatorname{exp}(D') < \operatorname{exp}(D)$.
	 \item\label{part:wreathproduct} $D \cong W_{t}$ for some $t \in \mathbb{N}$. 
	\end{enumerate}
\end{proposition}
\begin{proof} 
	\begin{enumerate}[(i)]
	\item From the definition of the fusion system of a block it is clear that $\mathcal{F}_b$ is a fusion subsystem of $\mathcal{F}_D(G)$. By Lemma \ref{lem:subsystem}, and Theorem \ref{thm:realisable} we are done.
	\item In this instance, the defect group $D$ is either abelian, or $D$ is maximal, so we are done by  Proposition \ref{prop:guarantee}(\ref{part:maximal}) and Corollary \ref{cor:abelian}. 
	\item Take $x \in D$ of maximal order $p^a$, and $\langle x \rangle$ fully-$\mathcal{F}_b$-centralized. Then if $x \in \Phi(C_D(x))$, we have $x = y^p z$, for some  $y \in C_D(x)$ and $z \in C_D(x)'$. Now as $z \in D'$, $z$ has order strictly less than $p^a$. Moreover as $y^p$ commutes with $x$, we have $x^{p^{a-1}}y^{-p^a} = z^{p^{a-1}}  = 1$. But this clearly implies that $x$ has order at most $p^{a-1}$, a contradiction. 
	\item The case $t=0$ is obvious. We show that $\operatorname{exp}(W_{t}') < \operatorname{exp}(W_{t})$ and apply part (\ref{part:lowerexponent}). We recall that the wreath product obeys the following short exact sequence:
		$$ 1 \rightarrow W_{t-1}^p \rightarrow W_{t} \rightarrow C_p \rightarrow 1.$$
		We immediately notice that $W_{t}' \leq  W_{t-1}^p$, and so $\operatorname{exp}(W_{t}') \leq \operatorname{exp}(W_{t-1})$. Write $p^a = \operatorname{exp}(W_{t-1})$. Given $x \in W_{t-1}$, of order $p^a$, it is clear the element $((x, e, \ldots, e), (1,\ldots ,p))$ has order $p^{a+1}$, so that $\operatorname{exp}(W_{t}) > \operatorname{exp}(W_{t-1}) \geq \operatorname{exp}(W_{t}')$. 
\end{enumerate}
\end{proof}
\subsection{Sporadic Groups}\label{sec:sporadic} Suppose that $S$ is a sporadic simple group, $E$ is a quasisimple extension of $S$ with $p'$-center, and $b$ is a block of $kE$. As $E$ has $p'$-center we have $|E|_p = |S|_p$. We see from \cite{ATLAS} that the only pair $(S, p)$ for which $p>5$ and $|S|_p>p^3$ is $(\mathrm{M},7)$. As such if $(S,p) \not= (\mathrm{M},7)$, $\mathcal{F}_b$ contains a fusion non-Schur element by Proposition \ref{prop:guarantee} (\ref{part:lessthanp3}). So we are only left with $(\mathrm{M},7)$. The Schur multiplier of $\mathrm{M}$ is trivial, so that $E = \mathrm{M}$. We see from \cite[Table 5.4]{WLA} that at $p=7$, all blocks of $\mathrm{M}$ are either of abelian or maximal defect, so we are done by Proposition \ref{prop:guarantee}(\ref{part:maximal}) and Corollary \ref{cor:abelian}. 
\subsection{Alternating Groups} Suppose that $n \in \mathbb{N}$, $E$ is a quasisimple extension of $A_n$ with $p'$-center, and $b$ is a block of $kE$ of defect $D$. If $n\leq 7$, then by the same reasoning as in subsection \ref{sec:sporadic}, $|E|_p = |A_n|_p \leq p$. As such, $D$ is abelian and we are done by Corollary \ref{cor:abelian}. So we can assume $n>7$, where we have $M(A_n)=2$. So then, either $E =A_n$ or $E =\tilde{A_n}$, the double cover of $A_n$. In the latter case, writing $\bar{b}$ for the image of $b$ in $kA_n$, $\bar{b}$ is either $0$, in which case $b$ is called faithful, or non-zero in which case $b$ is called non-faithful. 

Suppose that $b$ is a faithful block of $\tilde{A}_n$. Then by \cite[Proposition 2.4]{Gramain}, $D$ is a Sylow $p$-subgroup of $\tilde{A}_{m}$, for $m \leq n$. Note that the Sylow $p$-subgroups of $\tilde{A}_n, A_n$ and $S_n$ are isomorphic for $p>2$. If $b$ is non-faithful, the defect group $\bar{D}$ of $\bar{b}$ is isomorphic to $D$ by Proposition \ref{prop:usefulprop}(\ref{part:quotients}). Further, any block covering $\bar{b}$ of $S_n$ must also have defect group $D$. From \cite[Theorem 7.8]{Sambale} we know all blocks of $S_n$ have an associated \emph{weight} $\omega$, such that the defect group $D$ is isomorphic to a Sylow $p$-subgroup of $S_{p\omega}$. So in each case, (including when $E=A_n$), $D$ is isomorphic to a Sylow $p$-subgroup of $S_m$ for some $m \leq n$.

From \cite[Proposition 7.9]{Sambale}, each Sylow $p$-subgroup of $S_{m}$ is a direct product of groups of the form $W_{t}$. The exponent of a product of groups is the least common multiple of their exponents (which for $p$-powers is their maximum), so by Proposition \ref{prop:guarantee}(\ref{part:wreathproduct}), we see that $\mathcal{F}_b$ admits a fusion non-Schur element.
\section{Existence of fusion non-Schur elements for Lie-Chevalley groups}
Having dealt with the sporadic and alternating groups, we complete our proof by checking the criterion for quasisimple extensions of Lie-Chevalley groups. We briefly recount the structure of these groups.

Given a connected reductive linear algebraic group $\mathbf{G}$ defined over $\bar{\mathbb{F}}_r$, by \cite[Theorem 1.7]{MalleTesterman} there is an inclusion $i:\mathbf{G} \hookrightarrow \mathrm{GL}_n(\bar{\mathbb{F}}_r)$. Given $q=r^a$, we have an endomorphism $F_q$ of $\mathrm{GL}_n(\bar{\mathbb{F}}_r)$ sending $(x_{ij}) \mapsto (x_{ij}^q)$. A map $F:\mathbf{G}\rightarrow \mathbf{G}$ is called a Standard Frobenius endomorphism when $i \circ F = F_q \circ i$. More generally, a Frobenius endomorphism $F: \mathbf{G} \rightarrow \mathbf{G}$ is a map for which there exists $m \in \mathbb{N}$ with $F^m$ a standard Frobenius endomorphism. We refer to groups of the form $\mathbf{G}^F$ as finite groups of Lie type, with defining characteristic $r$. When $F$ is a standard Frobenius endomorphism with $i \circ F = F_q \circ i$, we say that $\mathbf{G}^F$ is defined over $\mathbb{F}_q$. The Lie-Chevalley groups are organised into the following families (stated in \cite{Wilson}).
\begin{enumerate}
	\item The classical groups: $\mathrm{PSL}_n(q), \mathrm{PSU}_n(q), \mathrm{PSp}_{2n}(q), \mathrm{P}\Omega_{2n+1}(q), \mathrm{P}\Omega_{2n}^+(q), \mathrm{P}\Omega_{2n}^-(q)$.
	\item The untwisted exceptional groups: $\mathrm{E}_6(q), \mathrm{E}_7(q), \mathrm{E}_8(q), \mathrm{F}_4(q), \mathrm{G}_2(q)$
	\item The twisted exceptional groups: $\mathrm{Sz}(2^{2n+1}), \hspace{0pt}^3\mathrm{D}_4(q),\hspace{0pt}^2\mathrm{E}_6(q), \hspace{0pt}^2\mathrm{F}_4(2^{2n+1})',\hspace{0pt}^2\mathrm{G}_2(3^{2n+1})$. 

\end{enumerate}
Each Lie-Chevalley group $S$ except $\mathrm{Sz}(2^{2n+1}),  \hspace{0pt}^2\mathrm{F}_4(2^{2n+1})'$ and $\hspace{0pt}^2\mathrm{G}_2(3^{2n+1})$ is realised as a central quotient of a finite group of Lie type, with a standard Frobenius endomorphism. Each of these groups is covered by a simply connected group of Lie type, which we write $G$, listed in \cite[Table 24.2]{MalleTesterman}. When $G$ is the universal perfect central extension of $S$, we say that $S$ has generic Schur multiplier, otherwise $S$ has exceptional Schur multiplier. The groups with exceptional Schur multiplier are listed in \cite[Table 24.3]{MalleTesterman}. Furthermore, we must distinguish between defining (where $r=p$) and non-defining characteristic. 

As such, we proceed as follows. First we verify the defining characteristic case, using the fact that all blocks have maximal or trivial defect. In non-defining characteristic, we deal with the Suzuki, Ree, and Tits cases, and the case where $G$ is not the universal cover. After this, we prove the remaining cases by appealing to a Jordan Decomposition result. 
\subsection{Defining Characteristic} 
Take a Lie-Chevalley group $S$, a quasisimple extension $E$ of $S$ with $p'$-center and a block $b$ of $kE$ of non-trivial defect. Suppose that $p=r$ and $p>5$. We first note that all Lie-Chevalley groups at $r>5$ have generic Schur multiplier, so that $E$ is a quotient of the simply connected cover $G$ of $S$; write $Z$ for the kernel. By \cite[Table 24.2]{MalleTesterman} we see $Z$ is a $p'$-group, and by Proposition \ref{prop:usefulprop}(\ref{part:quotients}) inflation from $E$ to $G$ preserves the defect group of $b$. From \cite[Theorem 1]{Humphreys}, the non-trivial defect groups of $p$-blocks of $G$ are maximal. As such, the defect group of $b$ is maximal in $E$, so that $\mathcal{F}_b$ contains a fusion non-Schur element, by Proposition \ref{prop:guarantee}(\ref{part:maximal}).

\subsection{Exceptional non-defining characteristic}
We complete the verification by checking the Lie-Chevalley groups in non-defining characteristic. We first deal with the Lie-Chevalley groups of exceptional Schur multiplier. 
\begin{proposition} Assume $p>5$. Suppose $S$ is a Lie-Chevalley group with exceptional Schur multiplier, $E$ a quasisimple extension of $S$ with $p'$-center and $b$ a block of $kE$ of non-trivial defect. Then $\mathcal{F}_b$ admits a fusion non-Schur element.
\end{proposition}
\begin{proof}
	Below is a table of the Lie-Chevalley groups with exceptional Schur multiplier, along with their Schur multiplier, and the maximal divisors of $|S|$ not containing $\{2,3,5\}$.
\begin{center}
	\resizebox{0.95\textwidth}{!}{%
	\begin{tabular}{|c|c|c|}
		\hline
			$S$ & $\mid M(S) \mid$  & $\mid S\mid_{\geq 7}$\\
	\hline
	\hline 
		$	\mathrm{SL}_2(4)$ &$ 2$ & $1$ \\
	\hline
		$	\mathrm{PSL}_2(9)$ &$ 6$ &$ 1$ \\
	\hline
		$	\mathrm{SL}_3(2)$ &$ 2$ &$ 7$ \\
	\hline

		$	\mathrm{PSL}_3(4)$ &$ 48$ &$ 7 $\\
	\hline
		$	\mathrm{SL}_4(2)$ &$ 2 $&$ 7 $\\
	\hline
		$	\mathrm{SU}_4(2)$ & $ 2$&$ 1$ \\ 
	\hline

	\end{tabular}
	\begin{tabular}{|c|c|c|}
		\hline
			$S$ &$ \mid M(S) \mid$  &$ \mid S\mid_{\geq 7}$\\
	\hline
	\hline 
			$\mathrm{PSU}_4(3)$ & $ 36$ & $7 $ \\
	\hline
			$\mathrm{PSU}_6(2)$ & $ 12$ &$ 7 \cdot 11$ \\
	\hline

			$\mathrm{Sp}_6(2)$ & $2$ & $7$ \\
	\hline
			$O_7(3)$ & $6$ & $7 \cdot 13 $ \\
	\hline
		$	SO_8^+(2) $& $4$ &  $7$ \\ 
	\hline
		$	\mathrm{Sz}(8) $&$4$ & $7 \cdot 13$ \\
	\hline

	\end{tabular}

	\begin{tabular}{|c|c|c|}
		\hline
			$S$ &$ \mid M(S) \mid$  &$ \mid S\mid_{\geq 7}$\\
	\hline

	\hline
			$ G_2(3)$ & $3$ & $7 \cdot 13$ \\
	\hline

			$ G_2(4) $ & $2$ & $7 \cdot 13$ \\
	\hline
			$ F_4(2)$ &$2$ & $7^2 \cdot 13 \cdot 17$\\
	\hline
			$ \hspace{0pt}^2E_6(2)$ & $12$  & $ 7^2 \cdot 11 \cdot 13 \cdot 17 \cdot 19 $\\ 
		  
	\hline 
		 	-&- &- \\
	\hline 		
			-& -& -\\
	\hline

	\end{tabular}
}
\end{center}
The data above, which can be found in \cite[Tables 24.1-24.3]{MalleTesterman}, show that in every case, $|S|_p \leq p^3$ and so $|E|_p \leq p^3$. Thus by Proposition \ref{prop:guarantee}(\ref{part:lessthanp3}) we are done.
\end{proof}

\subsection{The Suzuki, Ree and Tits Groups} The groups $\mathrm{Sz}(2^{2n+1})$  (at $n\geq 2$),  $\hspace{0pt}^2\mathrm{G}_2(3^{2n+1})$ and $\hspace{0pt}^2\mathrm{F}_4(2^{2n+1})$ have trivial Schur multiplier, so that in these cases there are no non-trivial quasisimple extensions. Further they have Weyl groups of order $2^3, 2^2 \cdot 3 $ and  $2^7 \cdot 3^2$ respectively (see \cite{Carter}). From \cite[Lemma 3.1 (2)]{FJL}, we see that at $p>5$, all Sylow $p$-subgroups of these groups are abelian, so by Corollary \ref{cor:abelian}, we are done. It remains to check $\mathrm{Sz}(8)$ and $\hspace{0pt}^2\mathrm{F}_4(2)'$. The first case was dealt with while considering the exceptional groups. In the latter case, the Tits group has order $2^{11} \cdot 3^3 \cdot 5^2 \cdot 13$, and so in particular has abelian Sylow $13$-subgroup. 
 
\subsection{Generic Lie-Chevalley groups in non-defining characteristic}
Suppose $S$ is a Lie-Chevalley group of defining characteristic $r$, defined over $\mathbb{F}_q$, with a quasisimple extension $E$ with $p'$-center and $G$ its universal perfect central extension.
Write $Z$ for the Sylow $p$-subgroup of $Z(G)$. We note that by Proposition \ref{prop:usefulprop}(\ref{part:quotients}), any block $\bar{b}$ of $kE$ inflates to a block $b$ of $kG/Z$, for which $\mathcal{F}_{\bar{b}}$ is a fusion subsystem of $\mathcal{F}_b$ on the same $p$-group. As such, it is enough to verify that the fusion systems of all blocks of $kG/Z$ admit fusion non-Schur elements. First, we recall some theory from \cite{Cabanes}.

Take $\mathbf{G}$ a connected reductive linear algebraic group defined over $\mathbb{F}_q$, with Frobenius endomorphism $F:\mathbf{G} \rightarrow \mathbf{G}$. We may define a group $\mathbf{G}^*$ dual to $\mathbf{G}$, admitting an endomorphism $F^*:\mathbf{G}^* \rightarrow \mathbf{G}^*$. To every semisimple $s \in (\mathbf{G}^*)^{F^*}$ there exists a subset $\mathcal{E}(\mathbf{G}^F, s) \subset \operatorname{Irr}(G^F)$. We write $\operatorname{Bl}(\mathbf{G}^F; s)$ for the blocks $b$ of $k\mathbf{G}^F$ so that $\operatorname{Irr}(b) \cap \mathcal{E}(\mathbf{G}^F, s) \not=\emptyset$. When $b \in \operatorname{Bl}(\mathbf{G}^F; 1)$, we call $b$ unipotent. If $\mathbf{G}$ is not connected, we say a block $b$ is unipotent when it covers a unipotent block of $(\mathbf{G}^\circ)^F$, the fixed points of the identity component of $\mathbf{G}$.
\begin{proposition}\cite[Theorem 7.11]{Cabanes}\label{prop:unipotent}
	Take $\mathbf{G}$ a connected reductive linear algebraic group, defined over $\mathbb{F}_q$ and $F: \mathbf{G} \rightarrow \mathbf{G}$ be the associated Frobenius endomorphism.  Assume $p \geq 7$ and $p \nmid q$. Take $b$ a unipotent block of $k\mathbf{G}^F$. Then $\mathcal{F}_b$ is realisable.
\end{proposition}

\begin{lemma}
	Take $\mathbf{G}$ a connected reductive linear algebraic group, defined over $\mathbb{F}_q$ and $F: \mathbf{G} \rightarrow \mathbf{G}$ be the associated Frobenius endomorphism.  Assume $p \geq 7$ and $p \nmid q$. Write $Z$ for the Sylow $p$-subgroup of $Z(\mathbf{G}^F$). Then, for every positive defect $p$-block $\bar{b}$ of $k\mathbf{G}^F/Z$, the block fusion system $\mathcal{F}_{\bar{b}}$ admits a fusion non-Schur element.
\end{lemma}
\begin{proof} 
	Take a block $b$ of $k\mathbf{G}^F$ of defect $D$ dominating $\bar{b}$. By Proposition \ref{prop:centralquotient} we have $\mathcal{F}_b/Z \cong \mathcal{F}_{\bar{b}}$. We first construct a realisable fusion system $\mathcal{E}$ on $D$, so that $\mathcal{F}_b$ is wrapped by $\mathcal{E}$. We then demonstrate that $C_\mathcal{E}(Z)/Z$ wraps $\mathcal{F}_b/Z$. We go on to justify that $C_\mathcal{E}(Z)/Z$ is realisable, and so we are done by Theorem \ref{thm:realisable}, and Lemma \ref{lem:automorph}.

	Take $s \in (\mathbf{G}^*)^{F^*}$ a semisimple $p'$-element such that $b \in \operatorname{Bl}(\mathbf{G}^F;s)$. From \cite[Theorem 1.4]{Enguehard}, there is a reductive group $\mathbf{G}(s)$, admitting a Frobenius endomorphism $\tilde{F}: \mathbf{G}(s) \rightarrow \mathbf{G}(s)$, and a fusion system preserving bijection between $\operatorname{Bl}(\mathbf{G}(s)^{\tilde{F}}; 1)$ and $\operatorname{Bl}(\mathbf{G}^F; s)$. This result requires that the Mackey formula holds for closed $F$-stable (resp. $F^*$-stable) subgroups of $\mathbf{G}$ (resp $\mathbf{G}^*$), which is true by \cite[Theorem 1.1]{Cotner}. Write $c$ for this corresponding unipotent block of $k\mathbf{G}(s)^{\tilde{F}}$. By assumption, $c$ covers a unipotent block $d$ of $k(\mathbf{G}(s)^\circ)^{\tilde{F}}$. Now, $(\mathbf{G}(s)^\circ)^{\tilde{F}}$ is a normal subgroup of $\mathbf{G}(s)^{\tilde{F}}$, which by \cite[Theorem 1.4(A), Proposition 1.2.4(a)]{Enguehard} has $[(\mathbf{G}(s))^{\tilde{F}}: (\mathbf{G}(s)^\circ)^{\tilde{F}}]$ coprime to $p$, so that $c$ and $d$ share the same defect group. Now, by Proposition \ref{prop:unipotent} we know $\mathcal{F}_d$ is realisable, and by Proposition \ref{prop:usefulprop}(\ref{part:normalimplication}) we know that $\mathcal{F}_d$ wraps $\mathcal{F}_b$. Write $\mathcal{E} :=\mathcal{F}_d$, and as $\mathcal{E}$ is realisable, choose a finite group $H$ containing $D$ as a Sylow $p$-subgroup for which $\mathcal{E} = \F_D(H)$. 

We now show that $C_\mathcal{E}(Z)/Z$ wraps $\F_b/Z$. Now, take $Z \leq R \leq D$, and a homomorphism $\bar\phi \in \Hom_{\F_b/Z}(R/Z, D/Z)$, induced by $\phi \in \Hom_{\F_b}(R,D)$ with $\phi(Z) = Z$.  Then we have $\psi \in \Hom_{\mathcal{E}}(R,D)$ and $\alpha \in \operatorname{Aut}(D)$ so that $\phi = \alpha \circ \psi$. As $\alpha \in \operatorname{Aut}(D)$ we know $\alpha^{-1}(Z) = \psi(Z)$ is central, and so fully-$\mathcal{E}$-centralized; therefore $\psi(Z)$ is receptive in $\mathcal{E}$. So by saturation of $\mathcal{E}$ we know $\psi|_{Z}:Z \rightarrow D$ extends to a morphism $N_{\psi|_{Z}} \rightarrow D$, where:
$$N_{\psi|_{Z}} := \{ d \in N_D(Z): \hspace{0pt}^{\psi|_{Z}} c_d \in \operatorname{Aut}_D(\psi(Z))\}=D.$$
So then, $\psi|_{Z}$ extends to an automorphism $\beta \in \operatorname{Aut}_\mathcal{E}(D)$. Write $\gamma = \alpha \circ \beta$ and $\theta = \beta^{-1} \circ \psi$. Then we have $\phi = \gamma \circ \theta$, where $\theta \in \Hom_{C_\mathcal{E}(Z)}(R,D)$. Moreover, $\gamma$ preserves $Z$, so that $\bar{\phi} = \bar{\gamma} \circ \bar{\theta}$, and we see $\mathcal{F}_b/Z$ is wrapped by $C_\mathcal{E}(Z)/Z$. 

Now, as we chose $H$ so that $\mathcal{E}= \F_D(H)$, by Proposition \ref{prop:centralizerblock} we know $C_\mathcal{E}(Z)\cong \mathcal{F}_D(C_H(Z))$. As $Z$ is normal in $C_H(Z)$, we know from \cite[Theorem 5.20]{Craven}, that $C_\mathcal{E}(Z)/Z \cong \F_{D/Z}(C_H(Z)/Z)$; that is, $C_\mathcal{E}(Z)/Z$ is realisable. Moreover, $D/Z$ is the defect group of $\bar{b}$, and so non-trivial. As such, we have $\mathrm{NS}(C_\mathcal{E}(Z)/Z)$. As $C_\mathcal{E}(Z)/Z$ wraps $\mathcal{F}_b/Z$, we know that $\mathcal{F}_{\bar{b}} \cong \mathcal{F}_b/Z$ admits a fusion non-Schur element.

\end{proof}

\subsection{Proof of existence of fusion non-Schur elements}
\begin{proof}[Proof of Theorem \ref{mainthm2}] The combined results of the last two sections verify that when $p>5$, if $E$ is a quasisimple group with $p'$-center, and $b$ is a block of $kE$ of non-trivial defect, that $\mathcal{F}_b$ admits a fusion non-Schur element. Then, Proposition \ref{prop:quasisimple} completes the proof.
\end{proof}
\bibliographystyle{alpha}
\bibliography{references}

@book {Wilson,
    AUTHOR = {Wilson, Robert A.},
     TITLE = {The finite simple groups},
    SERIES = {Graduate Texts in Mathematics},
    VOLUME = {251},
 PUBLISHER = {Springer-Verlag London, Ltd., London},
      YEAR = {2009},
     PAGES = {xvi+298},
      ISBN = {978-1-84800-987-5},
   MRCLASS = {20D05},
  MRNUMBER = {2562037},
MRREVIEWER = {Gernot\ Stroth},
       DOI = {10.1007/978-1-84800-988-2},
       URL = {https://doi.org/10.1007/978-1-84800-988-2},
}

@article {Carter,
    AUTHOR = {Carter, R. W.},
     TITLE = {Conjugacy classes in the {W}eyl group},
   JOURNAL = {Compositio Math.},
  FJOURNAL = {Compositio Mathematica},
    VOLUME = {25},
      YEAR = {1972},
     PAGES = {1--59},
      ISSN = {0010-437X,1570-5846},
   MRCLASS = {20H15 (17B20)},
  MRNUMBER = {318337},
MRREVIEWER = {Juri\ A.\ Bahturin},
}

@incollection {Linckelmann2,
    AUTHOR = {Linckelmann, Markus},
     TITLE = {Finite-dimensional algebras arising as blocks of finite group
              algebras},
 BOOKTITLE = {Representations of algebras},
    SERIES = {Contemp. Math.},
    VOLUME = {705},
     PAGES = {155--188},
 PUBLISHER = {Amer. Math. Soc., [Providence], RI},
      YEAR = {[2018] \copyright 2018},
      ISBN = {978-1-4704-3576-9},
   MRCLASS = {20C20 (16E40)},
  MRNUMBER = {3799494},
MRREVIEWER = {Constantin-Cosmin\ Todea},
       DOI = {10.1090/conm/705/14202},
       URL = {https://doi.org/10.1090/conm/705/14202},
}

@misc{LinckelmannTalk,
  author = {Linckelmann, M.},
  title = {On the {H}ochschild cohomology of finite group algebras},
  howpublished = "\url{https://www.newton.ac.uk/seminar/35694/}",
  year = {2022}, 
  note = "[Online; accessed October 2026]",
}

@book {Nagao,
    AUTHOR = {Nagao, Hirosi and Tsushima, Yukio},
     TITLE = {Representations of finite groups},
      NOTE = {Translated from the Japanese},
 PUBLISHER = {Academic Press, Inc., Boston, MA},
      YEAR = {1989},
     PAGES = {xviii+424},
      ISBN = {0-12-513660-9},
   MRCLASS = {20Cxx (20-01 20C20)},
  MRNUMBER = {998775},
MRREVIEWER = {Roderick\ Gow},
}

@misc{WLA,
      title={On $\mathrm{HH}^1$ of blocks of finite simple group algebras}, 
      author={Matthew Antrobus and Lleonard Rubio y Degrassi and Can Wen},
      year={2026},
      eprint={2611.12844},
      archivePrefix={arXiv},
      primaryClass={math.RT},
      url={https://arxiv.org/abs/2611.12844}, 
      note = {arXiv: \url{https://arxiv.org/abs/2611.12844}},
}

@misc{cotner,
      title={Base change and the {M}ackey formula in {D}eligne--{L}usztig theory}, 
      author={Sean Cotner},
      year={2026},
      eprint={2609.22331},
      archivePrefix={arXiv},
      primaryClass={math.RT},
      url={https://arxiv.org/abs/2609.22331}, 
	      note = {arXiv: \url{https://arxiv.org/abs/2609.22331}},
}

@misc{Enguehard,
      title={Towards a {J}ordan decomposition of blocks of finite reductive groups}, 
      author={Michel E. Enguehard},
      year={2013},
      eprint={1312.0106},
      archivePrefix={arXiv},
      primaryClass={math.GR},
      url={https://arxiv.org/abs/1312.0106}, 
      note = {arXiv: \url{https://arxiv.org/abs/1312.0106}},
}

@incollection {cabanes,
    AUTHOR = {Cabanes, Marc},
     TITLE = {Local methods for blocks of finite simple groups},
 BOOKTITLE = {Local representation theory and simple groups},
    SERIES = {EMS Ser. Lect. Math.},
     PAGES = {179--265},
 PUBLISHER = {Eur. Math. Soc., Z\"urich},
      YEAR = {2018},
      ISBN = {978-3-03719-185-9},
   MRCLASS = {20C20 (20D20)},
  MRNUMBER = {3821140},
MRREVIEWER = {Gabriel\ Navarro},
}

@article {Gramain,
    AUTHOR = {Gramain, Jean-Baptiste},
     TITLE = {The {I}saacs-{N}avarro conjecture for covering groups of the
              symmetric and alternating groups in odd characteristic},
   JOURNAL = {J. Algebraic Combin.},
  FJOURNAL = {Journal of Algebraic Combinatorics. An International Journal},
    VOLUME = {34},
      YEAR = {2011},
    NUMBER = {3},
     PAGES = {401--426},
      ISSN = {0925-9899,1572-9192},
   MRCLASS = {20C20 (05E10 20C15)},
  MRNUMBER = {2836368},
MRREVIEWER = {David\ Hemmer},
       DOI = {10.1007/s10801-011-0277-5},
       URL = {https://doi.org/10.1007/s10801-011-0277-5},
}

@article {Murphy2,
    AUTHOR = {Murphy, William},
     TITLE = {The {L}ie algebra structure of the degree one {H}ochschild
              cohomology of the blocks of the sporadic {M}athieu groups},
   JOURNAL = {J. Group Theory},
  FJOURNAL = {Journal of Group Theory},
    VOLUME = {26},
      YEAR = {2023},
    NUMBER = {1},
     PAGES = {161--191},
      ISSN = {1433-5883,1435-4446},
   MRCLASS = {20C34 (16E40 17B56 20C20)},
  MRNUMBER = {4529897},
MRREVIEWER = {Shigeo\ Koshitani},
       DOI = {10.1515/jgth-2021-0176},
       URL = {https://doi.org/10.1515/jgth-2021-0176},
}

@article {BKL2,
    AUTHOR = {Benson, David and Kessar, Radha and Linckelmann, Markus},
     TITLE = {Hochschild cohomology of symmetric groups and generating
              functions, {II}},
   JOURNAL = {Res. Math. Sci.},
  FJOURNAL = {Research in the Mathematical Sciences},
    VOLUME = {10},
      YEAR = {2023},
    NUMBER = {2},
     PAGES = {Paper No. 20, 5},
      ISSN = {2522-0144,2197-9847},
   MRCLASS = {20C20 (20J05)},
  MRNUMBER = {4588194},
MRREVIEWER = {Shigeo\ Koshitani},
       DOI = {10.1007/s40687-023-00382-2},
       URL = {https://doi.org/10.1007/s40687-023-00382-2},
}

@misc{Todea2,
      title={A reduction theorem for non-vanishing of {H}ochschild cohomology of block algebras and {H}appel's property}, 
      author={Patrick Serwene and Constantin-Cosmin Todea},
      year={2025},
      eprint={2503.05432},
      archivePrefix={arXiv},
      primaryClass={math.RT},
      url={https://arxiv.org/abs/2503.05432}, 
      note = {arXiv: \url{https://arxiv.org/abs/2503.05432}}, 
}

@article {Todea1,
    AUTHOR = {Todea, Constantin-Cosmin},
     TITLE = {Nontriviality of the first {H}ochschild cohomology of some
              block algebras of finite groups},
   JOURNAL = {J. Pure Appl. Algebra},
  FJOURNAL = {Journal of Pure and Applied Algebra},
    VOLUME = {227},
      YEAR = {2023},
    NUMBER = {2},
     PAGES = {Paper No. 107192, 10},
      ISSN = {0022-4049,1873-1376},
   MRCLASS = {16E40 (20C20)},
  MRNUMBER = {4457899},
MRREVIEWER = {Karin\ Erdmann},
       DOI = {10.1016/j.jpaa.2022.107192},
       URL = {https://doi.org/10.1016/j.jpaa.2022.107192},
}

@misc{Eisele,
      title={Block pro-fusion systems for profinite groups and blocks with infinite dihedral defect groups}, 
      author={Florian Eisele and Ricardo J. Franquiz Flores and John W. MacQuarrie},
      year={2025},
      eprint={2504.09286},
      archivePrefix={arXiv},
      primaryClass={math.RT},
      url={https://arxiv.org/abs/2504.09286}, 
      note = {arXiv: \url{https://arxiv.org/abs/2504.09286}}, 
}

@article {Humphreys,
    AUTHOR = {Humphreys, James E.},
     TITLE = {Defect groups for finite groups of {L}ie type},
   JOURNAL = {Math. Z.},
  FJOURNAL = {Mathematische Zeitschrift},
    VOLUME = {119},
      YEAR = {1971},
     PAGES = {149--152},
      ISSN = {0025-5874,1432-1823},
   MRCLASS = {20.70},
  MRNUMBER = {285623},
MRREVIEWER = {J.\ Sancho San Rom\'an},
       DOI = {10.1007/BF01109967},
       URL = {https://doi.org/10.1007/BF01109967},
}

@book {MalleTesterman,
    AUTHOR = {Malle, Gunter and Testerman, Donna},
     TITLE = {Linear algebraic groups and finite groups of {L}ie type},
    SERIES = {Cambridge Studies in Advanced Mathematics},
    VOLUME = {133},
 PUBLISHER = {Cambridge University Press, Cambridge},
      YEAR = {2011},
     PAGES = {xiv+309},
      ISBN = {978-1-107-00854-0},
   MRCLASS = {20G15 (20D06 20G05)},
  MRNUMBER = {2850737},
MRREVIEWER = {Timothy\ C.\ Burness},
       DOI = {10.1017/CBO9780511994777},
       URL = {https://doi.org/10.1017/CBO9780511994777},
}

@book {Sambale,
    AUTHOR = {Sambale, Benjamin},
     TITLE = {Blocks of finite groups and their invariants},
    SERIES = {Lecture Notes in Mathematics},
    VOLUME = {2127},
 PUBLISHER = {Springer, Cham},
      YEAR = {2014},
     PAGES = {xiv+243},
      ISBN = {978-3-319-12005-8; 978-3-319-12006-5},
   MRCLASS = {20C20 (20C15 20C40)},
  MRNUMBER = {3289382},
MRREVIEWER = {Charles\ W.\ Eaton},
       DOI = {10.1007/978-3-319-12006-5},
       URL = {https://doi.org/10.1007/978-3-319-12006-5},
}

@book {Craven,
    AUTHOR = {Craven, David A.},
     TITLE = {The theory of fusion systems},
    SERIES = {Cambridge Studies in Advanced Mathematics},
    VOLUME = {131},
      NOTE = {An algebraic approach},
 PUBLISHER = {Cambridge University Press, Cambridge},
      YEAR = {2011},
     PAGES = {xii+371},
      ISBN = {978-1-107-00596-9},
   MRCLASS = {20C20 (20D20 20D30)},
  MRNUMBER = {2808319},
MRREVIEWER = {Nadia\ P.\ Mazza},
       DOI = {10.1017/CBO9780511794506},
       URL = {https://doi.org/10.1017/CBO9780511794506},
}

@book {ATLAS,
    AUTHOR = {Conway, J. H. and Curtis, R. T. and Norton, S. P. and Parker,
              R. A. and Wilson, R. A.},
     TITLE = {{$\Bbb{ATLAS}$} of finite groups},
      NOTE = {Maximal subgroups and ordinary characters for simple groups,
              With computational assistance from J. G. Thackray},
 PUBLISHER = {Oxford University Press, Eynsham},
      YEAR = {1985},
     PAGES = {xxxiv+252},
      ISBN = {0-19-853199-0},
   MRCLASS = {20D05 (20-02)},
  MRNUMBER = {827219},
MRREVIEWER = {R.\ L.\ Griess},
}

@article{BKL,
	author = {Benson, Dave and Kessar, Radha and Linckelmann, Markus},
	doi = {10.2140/pjm.2021.313.1},
	fjournal = {Pacific Journal of Mathematics},
	issn = {0030-8730,1945-5844},
	journal = {Pacific J. Math.},
	mrclass = {16E40 (20C20)},
	mrnumber = {4313426},
	mrreviewer = {Karin\ Erdmann},
	number = {1},
	pages = {1--44},
	title = {On the {BV} structure of the {H}ochschild cohomology of finite group algebras},
	url = {https://doi.org/10.2140/pjm.2021.313.1},
	volume = {313},
	year = {2021}}

@article {KLR,
    AUTHOR = {Kessar, Radha and Linckelmann, Markus and Robinson, Geoffrey
              R.},
     TITLE = {Local control in fusion systems of {$p$}-blocks of finite
              groups},
   JOURNAL = {J. Algebra},
  FJOURNAL = {Journal of Algebra},
    VOLUME = {257},
      YEAR = {2002},
    NUMBER = {2},
     PAGES = {393--413},
      ISSN = {0021-8693,1090-266X},
   MRCLASS = {20C20},
  MRNUMBER = {1947328},
MRREVIEWER = {Ji\ Ping\ Zhang},
       DOI = {10.1016/S0021-8693(02)00517-3},
       URL = {https://doi.org/10.1016/S0021-8693(02)00517-3},
}

@article {Serwene,
    AUTHOR = {Serwene, Patrick},
     TITLE = {Block-exoticity of a family of exotic fusion systems},
   JOURNAL = {J. Group Theory},
  FJOURNAL = {Journal of Group Theory},
    VOLUME = {23},
      YEAR = {2020},
    NUMBER = {3},
     PAGES = {415--433},
      ISSN = {1433-5883,1435-4446},
   MRCLASS = {20D30 (20C20)},
  MRNUMBER = {4092937},
MRREVIEWER = {Charles\ W.\ Eaton},
       DOI = {10.1515/jgth-2019-0035},
       URL = {https://doi.org/10.1515/jgth-2019-0035},
}

@misc{BR,
      title={Outer derivations on blocks of group algebras}, 
      author={Matthew Antrobus and Benjamin Briggs and Lleonard Rubio y Degrassi},
      year={2026},
      eprint={2601.09602},
      archivePrefix={arXiv},
      primaryClass={math.RT},
      url={https://arxiv.org/abs/2601.09602}, 
      note = {arXiv: \url{https://arxiv.org/abs/2601.09602}}, 
}

@incollection{kulshammer,
	author = {K{\"u}lshammer, Burkhard},
	booktitle = {Representation theory of finite groups and finite-dimensional algebras ({B}ielefeld, 1991)},
	doi = {10.1007/978-3-0348-8658-1\_19},
	isbn = {3-7643-2604-2},
	mrclass = {16S34 (20C05 20C07 20C20)},
	mrnumber = {1112173},
	mrreviewer = {Karin\ Erdmann},
	pages = {425--442},
	publisher = {Birkh\"auser, Basel},
	series = {Progr. Math.},
	title = {Group-theoretical descriptions of ring-theoretical invariants of group algebras},
	url = {https://doi.org/10.1007/978-3-0348-8658-1_19},
	volume = {95},
	year = {1991}}

@book{linckelmann,
	author = {Linckelmann, Markus},
	isbn = {978-1-108-42590-2; 978-1-108-44180-3; 978-1-108-44190-2},
	mrclass = {20C20 (16-02 16D90 20-02)},
	mrnumber = {3821517},
	mrreviewer = {Burkhard\ K\"ulshammer},
	pages = {xi+510},
	publisher = {Cambridge University Press, Cambridge},
	series = {London Mathematical Society Student Texts},
	title = {The block theory of finite group algebras. {V}ol. {II}},
	volume = {92},
	year = {2018}}

@book {navarro,
    AUTHOR = {Navarro, G.},
     TITLE = {Characters and blocks of finite groups},
    SERIES = {London Mathematical Society Lecture Note Series},
    VOLUME = {250},
 PUBLISHER = {Cambridge University Press, Cambridge},
      YEAR = {1998},
     PAGES = {x+287},
      ISBN = {0-521-59513-4},
   MRCLASS = {20C20 (20-02 20C15)},
  MRNUMBER = {1632299},
MRREVIEWER = {Wolfgang\ Willems},
       DOI = {10.1017/CBO9780511526015},
       URL = {https://doi.org/10.1017/CBO9780511526015},
}

@book {AKO,
    AUTHOR = {Aschbacher, Michael and Kessar, Radha and Oliver, Bob},
     TITLE = {Fusion systems in algebra and topology},
    SERIES = {London Mathematical Society Lecture Note Series},
    VOLUME = {391},
 PUBLISHER = {Cambridge University Press, Cambridge},
      YEAR = {2011},
     PAGES = {vi+320},
      ISBN = {978-1-107-60100-0},
   MRCLASS = {20C20 (20D20 55P99)},
  MRNUMBER = {2848834},
MRREVIEWER = {Charles\ W.\ Eaton},
       DOI = {10.1017/CBO9781139003841},
       URL = {https://doi.org/10.1017/CBO9781139003841},
}

@article {KS,
    AUTHOR = {Kessar, Radha and Stancu, Radu},
     TITLE = {A reduction theorem for fusion systems of blocks},
   JOURNAL = {J. Algebra},
  FJOURNAL = {Journal of Algebra},
    VOLUME = {319},
      YEAR = {2008},
    NUMBER = {2},
     PAGES = {806--823},
      ISSN = {0021-8693,1090-266X},
   MRCLASS = {20C20 (20C15)},
  MRNUMBER = {2381811},
MRREVIEWER = {Tomoyuki\ Wada},
       DOI = {10.1016/j.jalgebra.2006.05.039},
       URL = {https://doi.org/10.1016/j.jalgebra.2006.05.039},
}

@article {Tradler,
    AUTHOR = {Tradler, Thomas},
     TITLE = {The {B}atalin-{V}ilkovisky algebra on {H}ochschild cohomology
              induced by infinity inner products},
   JOURNAL = {Ann. Inst. Fourier (Grenoble)},
  FJOURNAL = {Universit\'e{} de Grenoble. Annales de l'Institut Fourier},
    VOLUME = {58},
      YEAR = {2008},
    NUMBER = {7},
     PAGES = {2351--2379},
      ISSN = {0373-0956,1777-5310},
   MRCLASS = {16E40},
  MRNUMBER = {2498354},
MRREVIEWER = {Fran\c cois\ Goichot},
       DOI = {10.5802/aif.2417},
       URL = {https://doi.org/10.5802/aif.2417},
}

@article {FJL,
    AUTHOR = {Fleischmann, Peter and Janiszczak, Ingo and Lempken, Wolfgang},
     TITLE = {Finite groups have local non-{S}chur centralizers},
   JOURNAL = {Manuscripta Math.},
  FJOURNAL = {Manuscripta Mathematica},
    VOLUME = {80},
      YEAR = {1993},
    NUMBER = {2},
     PAGES = {213--224},
      ISSN = {0025-2611,1432-1785},
   MRCLASS = {20E34 (20C20 20D99)},
  MRNUMBER = {1233482},
MRREVIEWER = {Wolfgang\ Knapp},
       DOI = {10.1007/BF03026547},
       URL = {https://doi.org/10.1007/BF03026547},
}
\end{document}